\documentclass[11pt]{article}

\title{ Information on some recent applications of  umbral extensions to discrete mathematics}
\author{A.K.Kwa\'sniewski\\  
\\ Higher School of Mathematics and Applied Informatics\\
PL - 15-021 Bialystok , ul.Kamienna 17,  Poland
\\e-mail: kwandr@wp.pl}

\usepackage{amsmath,amsthm}

\chardef\bslash=`\\ 
\newtheorem{prop}{Proposition}[section]

\begin{document}
\maketitle\ \textit{to appear in  Review Bulletin of Calcutta
Mathematical Society Vol. 13 (2005)}
\begin{abstract}
At the first part of our communicate we show how specific umbral
extensions of the Stirling numbers of the second kind  result in
new type  of Dobinski-like formulas. In the second part among
others one recovers how and why  Morgan Ward solution of
uncountable family of  $\psi$- {\em difference calculus}
nonhomogeneous equations $\Delta_{\psi}f=\varphi$ in the form
$$f(x)=\sum_{n \geq
1}\frac{B_{n}}{n_{\psi}!}\varphi^{(n-1)}(x)+\int_{\psi}\varphi(x)
+p(x)$$ extends to $\psi$- Appell polynomials case automatically.
Illustrative specifications to $q$-calculus case and Fibonomial
calculus case are made explicit due to the usage of upside down
notation for objects of Extended Finite Operator Calculus .

\end{abstract}

MCS numbers:  05A40, 11B73, 81S99


\vspace{0.2cm} to be presented at the ISRAMA Congress, Calcuta -
India , December 2004

\vspace{0.5cm}
\subsection*{INTRODUCTION}

The upside down convenient notation for objects of Extended Finite
Operator Calculus  \textbf{(EFOC) } recently being developed and
promoted by the present author \cite{1,2,3,4} is to be used
throughout the whole exposition of both parts.\\

\vspace{2mm}

The \textbf{EFOC } is an implementation of operator approach of
Gian Carlo Rota to various by now formulations of extended of
umbral difference calculi. For abundant references - thousands  of
them - see \cite{5} and \cite{1,2,3,4}.

\vspace{2mm}

In \cite{1,2,3,4,5} and here also $\psi$ denotes  a number or
functions` sequence - sequence of functions of a parameter $q$.
$\psi$ constitutes the Morgan Ward [6] extension (see: after then
Viskov [7])  of $\langle\frac{1}{n!}\rangle_{n\geq 0}$ sequence to
quite arbitrary one (the so called - "admissible" in [8] and after
then in [1,2,3,4,5]).The specific choices are for example :
Fibonomialy-extended sequence $\langle\frac{1}{F_n!}\rangle_{n\geq
0}$ ($\langle F_n \rangle$ - Fibonacci sequence )  or just "the
usual" $\psi$-sequence $\langle\frac{1}{n!}\rangle_{n\geq 0}$ or
the famous Gauss $q$-extended $\langle\frac{1}{n_q!}\rangle_{n\geq
0}$ admissible sequence of extended umbral operator calculus,
where $ n_q=\frac{1-q^n}{1-q}$ and $n_q!=n_q(n-1)_q! , 0_q!=1$ -
see more below. With such type extension we frequently may in a
"$\psi$-mnemonic" way repeat reasoning and this what was done by
Rota (see at first \textbf{(FOC) } of Rota in [9]).

\vspace{2mm}

The simplicity of the first steps to be done while identifying
general properties of  $\psi$-extended objects  consists in
writing objects of these extensions in mnemonic convenient
\textbf{upside down notation} \cite{1,2,3,4} , \cite{10} which is
here down introduced as follows:

\begin{equation}\label{eq1}
\frac {\psi_{(n-1)}}{\psi_n}\equiv n_\psi,
n_\psi!=n_\psi(n-1)_\psi!, n>0 ,   x_{\psi}\equiv \frac
{\psi{(x-1)}}{\psi(x)} ,
\end{equation}

\begin{equation}\label{eq2}
x_{\psi}^{\underline{k}}=x_{\psi}(x-1)_\psi(x-2)_{\psi}...(x-k+1)_{\psi}
\end{equation}

\begin{equation}\label{eq3}
x_{\psi}(x-1)_{\psi}...(x-k+1)_{\psi}=
\frac{\psi(x-1)\psi(x-2)...\psi(x-k)} {\psi(x)
\psi(x-1)...\psi(x-k +1)} .
\end{equation}

If one writes the above in the form $x_{\psi} \equiv \frac
{\psi{(x-1)}}{\psi(x)}\equiv \Phi(x)\equiv\Phi_x\equiv x_{\Phi}$ ,
one sees that the name upside down notation is legitimate.

 You may consult \cite{1,2,3,4,5} and \cite{10} for further
development and use of this notation .

\vspace{2mm}

In the \textbf{first} part of our communicate umbral extensions of
the Stirling numbers of the both kinds are considered  and the
resulting new type  of Dobinski-like formulas are discovered [10].
These extensions naturally encompass the well known $q$-extensions
. The fact that $q$-extended Stirling numbers giving rise to the
umbral $q$-extended Dobinski formula interpreted as the average of
powers of random variable $X_q$ with the $q$-Poisson distribution
and are equivalent by re-scaling with the other Comtet -like
$q$-extended Stirling numbers [10] -  singles out the
$q$-extensions itself which appear to be a kind of bifurcation
point in the domain of umbral extensions . The further consecutive
umbral extensions of Carlitz-Gould $q$-Stirling numbers are
therefore realized in [10] in a two-fold way.

\vspace{2mm}

In the \textbf{second } part of our communicate one displays
\cite{11} how and  why Morgan Ward solution \cite{6} of $\psi$-
{\em difference calculus} nonhomogeneous equation
$\Delta_{\psi}f=\varphi$ in the form
$$f(x)=\sum_{n \geq
1}\frac{B_{n}}{n_{\psi}!}\varphi^{(n-1)}(x)+\int_{\psi}\varphi(x)
+p(x)$$ recently proposed by the present author \cite{12}- extends
here now to $\psi$- Appell polynomials case - almost {\em
automatically}. Illustrative specifications to $q$-calculus case
and Fibonomial calculus case \cite{2,13,14} were already made
explicit in \cite{12} exactly due to the of upside down notation
for objects of the  \textbf{EFOC }.

\vspace{2mm}

\subsection*{The First Part - Stirling numbers extensions and Dobinski-like formulae}

In this part of our communicate we follow [10] and we refer Reader
for further details to consult [10]. At the start let us recall
that the \textbf{two} standard \cite{15,16,17} $q$-extensions
Stirling numbers of the second kind  might be defined as follows:

\begin{equation}\label{eq4}
x_q^n=\sum_{k=0}^{n}\Big\{{n \atop k}\Big\}_q  x_q^{\underline k},
\end{equation}

where $ x_q=\frac{1-q^x}{1-q}$ and $x_q^{\underline
k}=x_q(x-1)_q...(x-k+1)_q $ , which corresponds to the $\psi$
sequence choice in the $q$-Gauss form
$\langle\frac{1}{n_q!}\rangle_{n\geq 0}$

\vspace{2mm}

and $q^\sim$-Stirling numbers

\begin{equation}\label{eq5}
x^n=\sum_{k=0}^{n}\Big\{{n \atop k}\Big\}^\sim_q \chi_{\underline
k}(x)
\end{equation}

where $\chi_{\underline k}(x)= x(x-1_q)(x-2_q)...(x-[k-1]_q)$

For these two classical by now $q$-extensions of  Stirling numbers
of the second kind  - the "$q$-standard" recurrences hold
respectively:

$$ \Big\{{{n+1}\atop k}\Big\}_q =
\sum_{l=0}^{n}\binom{n}{l}_q q^l\Big\{{l\atop {k-1}}\Big\}_q ;
n\geq 0 , k\geq 1,$$

$$ \Big\{{{n+1}\atop k}\Big\}^\sim_q =
\sum_{l=0}^{n}\binom{n}{l}_q q^{l-k+1}\Big\{{l\atop
{k-1}}\Big\}^\sim_q  ;  n\geq 0 , k\geq 1.$$

From the above it follows immediately that corresponding
$q$-extensions of  $B_n$ Bell numbers satisfy respective
recurrences:

$$ B_q(n+1) =
\sum_{l=0}^{n}\binom{n}{l}_q q^l B_q(l) ; n\geq 0 ,$$

$$ B^\sim_q(n+1) =
\sum_{l=0}^{n}\binom{n}{l}_q q^{l-k+1} \overline {B} ^\sim_q(l)  ;
n\geq 0 $$ where
$$ \overline {B}^\sim_q(l)=
\sum_{k=0}^{l} q^k\Big\{{l\atop k}\Big\}^\sim_q .
$$

\textbf{Note} that both definitions via (4) and (5) equations
consequently correspond to different $q$-counting   \cite{16} .\\

For applications to coherent state phisics see [17]  and
references therein.\\

 With any other choice out of countless choices of the
$\psi$ sequence the equation (5) becomes the definition of
$\psi^\sim$- Stirling (vide "Fibonomial-Stirling") numbers of the
second kind$\Big\{{n \atop k}\Big\}^\sim_\psi$  and then
$\psi^\sim$-Bell numbers $B^\sim_n(\psi)$ are defined as sums as
usual - where now $\chi_{\underline k}(x)$ in (5) is to be
replaced by $\psi_{\underline k}(x)=
x(x-1_\psi)(x-2_\psi)...(x-[k-1]_\psi)$. These $\psi^\sim$-
Stirling numbers of the second kind recognized for  $q$  case
properly as Comtet numbers in Wagner`s terminology \cite{16}
satisfy familiar recursion and are given by familiar formulas to
be presented soon.\\

The extension of definition (4) of the $q$- Stirling numbers of
the second kind beyond this $q$-case i.e. beyond the
$\psi=\langle\frac{1}{n_q!}\rangle_{n\geq 0}$ choice is not that
mnemonic and the " behavior " of the naturally  expected recursion
under extension  comprises quite a surprise (see: appendix 2.2 in
[10]).\\

Therefore  further consecutive umbral extension of Carlitz-Gould
$q$-Stirling numbers   $\Big\{{n\atop k}\Big\}_q $  and
$\Big\{{n\atop k}\Big\}^\sim_q $   is realized two-fold way [10].

The first "easy way" consists in almost mnemonic sometimes
replacement of $q$ subscript by $\psi$ because we learn from
\cite{16} that the equation (5) defines as a matter of fact the
specific case of Comtet numbers \cite{16} i.e. we define
$\psi$-extended Comtet-Stirling numbers of the second kind as
follows:

\begin{equation}\label{eq6}
x^n=\sum_{k=0}^{n}\Big\{{n \atop k}\Big\}^\sim_{\psi}
\psi_{\underline k}(x)
\end{equation}

where $\psi_{\underline k}(x)=
x(x-1_{\psi})(x-2_{\psi})...(x-[k-1]_{\psi}).$

As a consequence we have "for granted" :

\begin{equation}\label{eq7}
\Big\{{{n+1}\atop k}\Big\}^\sim_{\psi} = \Big\{{n\atop
{k-1}}\Big\}^\sim_{\psi} + k_{\psi}\Big\{{n\atop
k}\Big\}^\sim_{\psi} ;\quad n\geq 0 , k\geq 1 ;
\end{equation}

where \quad $\Big\{{n\atop 0}\Big\}^\sim_{\psi}=
\delta_{n,0},\quad\Big\{{n\atop k}\Big\}^\sim_{\psi}=0 ,\quad k>n
;\quad $\quad and the easy derivable  recurrence relations for
ordinary generating function now read

\vspace{2mm}

\begin{equation}\label{eq8}
G^\sim_{k_{\psi}}(x)=\frac{x}{1-k_{\psi}}G^\sim_{k_{\psi}-1}(x) ,
\quad k\geq 1
\end{equation}

where
$$G^\sim_{k_{\psi}}(x)= \sum_{n\geq 0}\Big\{{n\atop
k}\Big\}^\sim_{\psi}x^n     ,\quad k\geq 1 $$

\vspace{2mm}

from where as a consequence we arrive to what follows:
\vspace{5mm}

\begin{equation}\label{eq9}
G^\sim_{k_{\psi}}(x)=\frac{x^k}{(1-1_{\psi}x)(1-2_{\psi}x)...(1-k_{\psi}x)}
\quad, \quad k\geq 0
\end{equation}

\vspace{2mm}

and adapting reasoning from  \cite{18} we derive the following
explicit formula

\begin{equation}\label{eq10}
\Big\{{n\atop k}\Big\}^\sim_{\psi} = \frac
{r_{\psi}^n}{k_{\psi}!}\sum_{r=1}^{k}(-1)^{k-r}\binom{k_{\psi}}{r_{\psi}};
\quad n,k\geq 0 .
\end{equation}

\vspace{2mm}

Expanding the right hand side of (9) results in another explicit
formula for these $\psi$-case Comtet numbers i.e. we have

\vspace{2mm}

\begin{equation}\label{eq11}
\Big\{{n\atop k}\Big\}^\sim_{\psi} = \sum_{1\leq i_1 \leq
i_2\leq...\leq i_{n-k}\leq
k}(i_1)_{\psi}(i_2)_{\psi}...(i_{n-k})_{\psi}; \quad n,k\geq 0 .
\end{equation}

\vspace{2mm}

or equivalently  (compare with [16])

 \vspace{2mm}
\begin{equation}\label{eq12}
\Big\{{n\atop k}\Big\}^\sim_{\psi}=\sum_{d_1+ d_2+...+d_k =
n-k,\quad d_i\geq
0}1_{\psi}^{d_1}2_{\psi}^{d_2}...k_{\psi}^{d_k};\quad n,k\geq 0 .
\end{equation}

\vspace{2mm}

With help of $\psi^\sim$-\textbf{Stirling numbers} of the second
kind being defined equivalently by (6) , (7), (11)  or (12) we
define now  $\psi^\sim$-\textbf{Bell numbers} in a standard way
$$ B^\sim_n(\psi)=\sum_{k=0}^{n} \Big\{{n\atop
k}\Big\}^\sim_{\psi} ,\qquad n\geq 0 .$$

The recurrence  for  $B^\sim_n(\psi)$ is already quite involved
and complicated for the  $q$-extension case (see: the first
section in [10])- and no acceptable readable form of recurrence
for the $\psi$-extension case is known to us.

Nevertheless after adapting the  the corresponding Wilf`s
reasoning from \cite{18} we derive  for two variable ordinary
generating function for $\Big\{{n\atop k}\Big\}^\sim_{\psi}$
Stirling numbers of the second kind and the $\psi$-exponential
generating function for $ B^\sim_n(\psi)$  Bell numbers the
following  formulae

\begin{equation}\label{eq13}
C^\sim_{\psi}(x,y) = \sum_{n\geq 0} A^\sim_n (\psi,y)x^n ,
\end{equation}
where the $\psi$- exponential-like polynomials $ A^\sim_n
(\psi,y)$
$$ A^\sim_n (\psi,y)=\sum_{k=0}^{n} \Big\{{n\atop
k}\Big\}^\sim_{\psi}y^k$$ do satisfy the recurrence

$$ A^\sim_n (\psi,y)=  [y(1+\partial_{\psi}]A^\sim_{n-1}(\psi,y) \qquad n\geq 1 ,$$

hence
$$ A^\sim_n (\psi,y)=  [y(1+\partial_{\psi}]^n 1,\quad \qquad n\geq 0 ,$$
where  the linear operator $\partial_{\psi}$ acting on the algebra
of formal power series is being called (see: [1,2,3,4,5] and
references therein) the "$\psi$-derivative" and $\partial_{\psi}
y^n = n_{\psi}y^{n-1}.$

\vspace{2mm}

The $\psi$-exponential generating function $ F^{\sim}_{\psi}(x)=
\sum_{n\geq 0}B^\sim_n(\psi)\frac{x^n}{n_{\psi}!} $

for $ F^{\sim}_n(\psi)$ Bell numbers - again - after cautious
adaptation of the method from the Wilf`s generatingfunctionology
book  [18] we get only a little bit involved formula

\begin{equation}\label{eq14}
B^{\sim}_{\psi}(x)= \sum_{r\geq 0}\epsilon(\psi,r)
\frac{e_{\psi}[r_{\psi}x]}{r_{\psi}!}
\end{equation}

where (see: [6,7,1,2,3,4])
$$e_{\psi}(x) =
\sum_{n\geq 0}\frac{x^n}{n_{\psi}!}$$

while
\begin{equation}\label{eq15}
\epsilon(\psi,r)=\sum_{k=r}^{\infty} \frac{(-1)^{k-r}}{(k_{\psi}-
r_{\psi})!}
\end{equation}

and  for  $\psi$-extension the Dobinski like formula here now
reads

\begin{equation}\label{eq16}
B^{\sim}_n (\psi)= \sum_{r\geq 0}\epsilon(\psi,r)
\frac{r_{\psi}^n}{r_{\psi}!}.
\end{equation}

In the case of  Gauss $q$-extended  choice of
$\langle\frac{1}{n_q!}\rangle_{n\geq 0}$ admissible sequence of
extended umbral operator calculus  we have then

\begin{equation}\label{eq17}
\epsilon(q,r)=\sum_{k=r}^{\infty}
\frac{(-1)^{k-r}}{(k-r)_q!}q^{-\binom {r}{2}}
\end{equation}
and the new $q^\sim$-Dobinski  formula is given by

\begin{equation}\label{eq18}
B^{\sim}_n (\psi)= \sum_{r\geq 0}\epsilon(\psi,r)
\frac{r_{\psi}^n}{r_{\psi}!}.
\end{equation}
which for $ q=1$ becomes the Dobinski formula from 1887 [19].\\

As for the problem of how eventually one might  interpret the
$\psi^\sim$-Dobinski formulae (16) and (18) in the Rota-like way
see: [10].\\

In [10] you  find also details in support of  our conviction that
$q$-extensions seem  appear as a kind of bifurcation point in the
domain of umbral extensions.\\

For the discussion of the other way to arrive eventually to
another type of Dobinski formula  - via an attempt to
$\psi$-extend the equation (4) - we again refer the Reader to [10]
as this is quite a longer story with a surprise.

\vspace{2mm}

The parallel treatment  of the Comtet $ \left[ {{\begin{array}{*{20}c} {n}\\
{k}\end{array}} } \right]^\sim_\psi$ Stirling numbers of the first
kind is now not difficult.

Namely we define  as in [10] the $\psi^\sim$ Stirling numbers of
the first kind as the coefficients in the following expansion

\begin{equation}\label{eq19}
\psi_{\underline k}(x)=\sum_{r=0}^{k}\left[ {{\begin{array}{*{20}c} {k}\\
{r}\end{array}} } \right]^\sim_\psi x^r
\end{equation}

where -  recall $\psi_{\underline k}(x)=
x(x-1_{\psi})(x-2_{\psi})...(x-[k-1]_{\psi});\quad$ From the above
we then get

\begin{equation}\label{eq20}
\sum_{r=0}^{k}\left[ {{\begin{array}{*{20}c} {k}\\
{r}\end{array}} } \right]^\sim_\psi \Big\{{r\atop
l}\Big\}^\sim_{\psi}= \delta_{k,l}.
\end{equation}

Another (expected Whithney numbers of the first kind) $\psi^c$-
Stirling numbers of the first kind [10]  are defined as follows:

\begin{equation}\label{eq21}
\psi_{\overline k}(x)=\sum_{r=0}^{k}\left[ {{\begin{array}{*{20}c} {k}\\
{r}\end{array}} } \right]^c_\psi x^r
\end{equation}

where -  now  $\psi_{\overline k}(x)=
x(x+1_{\psi})(x+2_{\psi})...(x+[k-1]_{\psi});\quad$ More on that -
see [10].

\subsection*{The Second  Part -$\psi$- Appell polynomials` solutions of
the $Q(\partial_{\psi})$- difference  nonhomogeneous equation}

The second part is based on [11,12] to which we refer for more
details.

\vspace{2mm}

At first  recall \cite{3,1} the simple fact to be used in what
follows.

\begin{prop}
$Q(\partial_{\psi})$ is a $\psi$- delta operator iff there exists
invertible $S\in \Sigma _{\psi}$ such that
$Q(\partial_{\psi})=\partial_{\psi}S$.
\end{prop}

Formally: "$S=Q/\partial_{\psi}$" or "$S^{-1}=\partial_{\psi}/Q$".
In the sequel we use this abbreviation $Q(\partial_{\psi})\equiv
Q$.

$\psi$- {\em Appell} or generalized {\em Appell polynomials}
$\left\{A_{n}(x)\right\}_{n \geq 0}$ are defined according to
\begin{equation} \label{psi-A-1}
\partial_{\psi}A_{n}(x)=n_{\psi}A_{n-1}(x)
\end{equation}
and they naturally do satisfy the $\psi$ - Sheffer-Appell identity
\cite{3,1}
\begin{equation}\label{psi-A-2}
A_{n}(x+_{\psi}y)=\sum_{s=0}^{n}\binom{n}{s}
_{\psi}A_{s}(y)x^{n-s}.
\end{equation}
$\psi$-{\em Appell} or generalized {\em Appell polynomials}
$\left\{A_{n}(x)\right\}_{n \geq 0}$ are equivalently
characterized via their $\psi$- exponential generating function
\begin{equation}\label{psi-A-3}
\sum_{n \geq
0}z^{n}\frac{A_{n}(x)}{n_{\psi}!}=A(z)\exp_{\psi}\left\{xz\right\},
\end{equation}
where $A(z)$ is a formal series with constant term different from
zero - here normalized to one.

The $\psi$- exponential function of $\psi$-Appell-Ward numbers
$A_{n}=A_{n}(0)$ is
\begin{equation} \label{psi-A-4}
\sum_{n \geq 0}z^{n}\frac{A_{n}}{n_{\psi}!}=A(z).
\end{equation}
Naturally $\psi$- {\em Appell} $\left\{A_{n}(x)\right\}_{n \geq
0}$ satisfy the $\psi$- {\em difference} equation
\begin{equation}\label{psi-A-5}
QA_{n}(x)=n_{\psi}x^{n-1};\;\;\;'n\geq 0,
\end{equation}
because
$QA_{n}(x)=QS^{-1}x^{n}=Q(\partial_{\psi}/Q)x^{n}=\partial_{\psi}x^{n}=n_{\psi}x^{n-1}\;
;n \geq 0$. Therefore they play the same role in
$Q(\partial_{\psi})$- {\em difference} calculus as Bernoulli
polynomials do in standard difference calculus or $\psi$-{\em
Bernoulli-Ward polynomials} (see Theorem 16.1 in \cite{1} and
consult also \cite{12}) in $\psi$-{\em difference} calculus due to
the following: The central problem of the $Q(\partial_{\psi})$ -
{\em difference calculus} is:
$$ Q(\partial_{\psi})f=\varphi\;\;\;\;\;\;\;\;\;\;\varphi=?,$$
where $f, \varphi$ - are for example formal series or polynomials.

The idea of finding solutions is the $\psi$-{\em Finite Operator
Calculus} \cite{1,3,4,2} standard. As one knows  \cite{1,3}) any
$\psi$- delta operator $Q$ is of the form

$$Q(\partial_{\psi})=\partial_{\psi}S $$

where$ S\in \Sigma_{\psi}.$ Let $Q(\partial_{\psi})=\sum_{k \geq
1}\frac{q_{k}}{k_{\psi}!}\partial_{\psi}^{k},\;\;q_{1}\neq 0$.
Consider then $Q(\partial_{\psi}=\partial_{\psi}S)$ with
$S=\sum_{k \geq
0}\frac{q_{k+1}}{(k+1)_{\psi}!}\partial_{\psi}^{k}\equiv \sum_{k
\geq
0}\frac{s_{k}}{k_{\psi}!}\partial_{\psi}^{k};\;\;\;s_{0}=q_{1}\neq
0$.

We thus have for $S^{-1}\equiv\hat{A}$ - call it: $\psi$- {\em
Appell operator} - the obvious expression

$$\hat{A}\equiv S^{-1}=\frac{\partial_{\psi}}{Q_{\psi}}=\sum_{n\geq 0}\frac{A_{n}}{n_{\psi}!}\partial_{\psi}^{n}.$$

Now multiply the equation \quad $Q(\partial_{\psi})f=\varphi
$\quad  by $\hat{A}\equiv \sum_{n\geq
0}\frac{A_{n}}{n_{\psi}!}\partial_{\psi}^{n}$ thus getting

\begin{equation}\label{psi_A_6}
\partial_{\psi}f=\sum_{n\geq
0}\frac{A_{n}}{n_{\psi}!}\varphi^{(n)},\;\;\;\varphi^{(n)}=\partial_{\psi}\varphi^{(n-1)}.
\end{equation}
The solution then reads:

\begin{equation}\label{psi-A-7}
f(x)=\sum_{n \geq
1}\frac{A_{n}}{n_{\psi}!}\varphi^{(n-1)}(x)+\int_{\psi}\varphi(x)+p(x),
\end{equation}

where $p$ is "$Q(\partial_{\psi})$- periodic" i.e.
$Q(\partial_{\psi})p=0$. Compare with  \cite{12}

for "$+_{\psi}1$- periodic" i.e. $p(x+_{\psi}1)=p(x)$ i.e.
$\Delta_{\psi}p=0$. Here
 the relevant  $\psi$ - integration $\int_{\psi}\varphi(x)$ is defined as in
\cite{1}. We recall it in brief. Let us introduce the following
representation for $\partial_{\psi}$ "difference-ization"

$$ \partial_{\psi}=\hat{n}_{\psi}\partial_{0}\; ;\;\;\;\hat{n}_{\psi}x^{n-1}=n_{\psi}x^{n-1};\;\;n\geq 1,$$

where $\partial_{0}x^{n}=x^{n-1}$ i.e. $\partial_{0}$ is the $q=0$
Jackson derivative. $\partial_{0}$ is identical with divided
difference operator. Then we define the linear mapping
$\int_{\psi}$ accordingly:

$$ \int_{\psi}x^{n}=\left( \hat{x}\frac{1}{\hat{n}_{\psi}}\right)x^{n}=\frac{1}{(n+1)_{\psi}}x^{n+1};\;\;\;n\geq 0$$
where of course  $\partial_{\psi}\circ \int_{\psi}=id$.

\section{Examples}
\renewcommand{\labelenumi}{(\alph{enumi})}
\begin{enumerate}

\item The case of $\psi$- {\em Bernoulli-Ward} polynomials and
$\Delta_{\psi}$- {\em difference calculus} was considered in
detail in \cite{12} following \cite{6}. \item Specification of (a)
to the Gauss and Heine originating $q$-umbral calculus case
\cite{6,3,4,5,2} was already presented in \cite{12}.

\item Specification of (a) to the Lucas originating FFOC - case
was also presented in \cite{12} (here: FFOC={\bf F}ibonomial {\bf
F}inite {\bf O}perator {\bf C}alculus), see example 2.1 in
\cite{2}). Recall: the {\em Fibonomial coefficients} -already
known in 19-th century to Lucas [20]-  are defined "binomially" as

$$\binom{n}{k}_{F}=\frac{F_{n}!}{F_{k}!F_{n-k}!}=\binom{n}{n-k}_{F},$$

where ($F_{n}$- {\em Fibonacci numbers} in up-side down notation:
$n_{F}\equiv F_{n}\neq
0$,\\$n_{F}!=n_{F}(n-1)_{F}(n-2)_{F}(n-3)_{F}\ldots
2_{F}1_{F};\;\;0_{F}!=1$;\\$n_{F}^{\underline{k}}=n_{F}(n-1)_{F}\ldots
(n-k+1)_{F};\;\;\;\binom{n}{k}_{F}\equiv
\frac{n^{\underline{k}}_{F}}{k_{F}!}$.

\vspace{2mm}

We shall call the corresponding linear difference operator
$\partial_{F};\;\;\partial_{F}x^{n}=n_{F}x^{n-1};\;\;n\geq 0$ the
$F$-derivative. Then in conformity with \cite{6} and with notation
as in \cite{1}-\cite{13,14} one has:

$$E^{a}(\partial_{F})=\sum_{n\geq 0}\frac{a^{n}}{n_{F}!}\partial_{F}^{n}$$

for the corresponding [3,1,4] generalized translation operator
$E^{a}(\partial_{F})$. The $\psi$- integration for the moment is
still not explored $F$- integration and we arrive at the
\textbf{$F$-
Bernoulli} polynomials unknown till now.\\

{\bf Note:} recently a combinatorial interpretation of Fibonomial
coefficient has been found \cite{14} by the present author.

\item The other examples of $Q(\partial_{\psi})$- {\em difference
calculus} - expected naturally to be of primary importance in
applications are provided by the possible use of such $\psi$-
Appell polynomials as:

\begin{itemize}
\item

$\psi$-Hermite polynomials $\left\{H_{n,\psi}\right\}_{n\geq 0}$:

$$ H_{n,\psi}(x)=\left[\sum_{k\geq 0}\left(-\frac{1}{2}\right)^{k}
\frac{\partial_{\psi}^{2k}}{k_{\psi}!}\right]x^{n}\;\;\;n\geq 0;$$

\item

$\psi$ - Laguerre polynomials $\left\{ L_{n,\psi}\right\}_{n\geq
0}$ \cite{3}:

\begin{multline*}
L_{n,q}(x) = \frac{{n_{q}} }{{n}}\hat {x}_{\psi}
\left[{\frac{{1}}{{\partial _{\psi} - 1}}}\right]^{-n}x^{n-1} =
\frac{{n_{\psi}} }{{n}}\hat {x}_{\psi} \left( {\partial _{\psi} -
1} \right)^{n}x^{n-1} = \\ = \frac{{n_{\psi}} }{{n}}\hat
{x}_{\psi} \sum\limits_{k = 1}^{n} \left( { - 1} \right)^{k}\left(
{{\begin{array}{*{20}c}
 {n} \hfill \\
 {k} \hfill \\
\end{array}} } \right) \partial _{\psi} ^{n - k}x^{n - 1} =\\
=\frac{{n_{\psi}} }{{n}}\sum\limits_{k = 1}^{n} {} \left( { - 1}
\right)^{k}\left( {{\begin{array}{*{20}c}
 {n} \hfill \\
 {k} \hfill \\
\end{array}} } \right) \left( {n - 1} \right)_{\psi} ^{\underline {n - k}
}\frac{{k}}{{k_{\psi}} }x^{k} .
\end{multline*}

For $q=1$ in $q$-extended case one recovers the known formula :

\begin{center}
$L_{n,q=1}(x) = \sum\limits_{k = 1}^{n} {} \left( { - 1}
\right)^{k}\frac{{n_{q} !}}{{k_{q} !}}\left(
{{\begin{array}{*{20}c}
 {n - 1} \hfill \\
 {k - 1} \hfill \\
\end{array}} } \right)x^{k}$.
\end{center}

\end{itemize}
\end{enumerate}

\end{document}